\documentclass[a4paper,12pt]{amsart}

\usepackage[dvipdfmx]{graphicx}
\usepackage[top=30truemm,bottom=30truemm,left=25truemm,right=22truemm]{geometry}
\usepackage{comment}
\usepackage{bm}
\usepackage{xcolor}
\usepackage{amsmath,amssymb}
\usepackage{latexsym}
\usepackage{amsthm}
\numberwithin{equation}{section}

\usepackage{amsaddr}

\usepackage{amsrefs}

\theoremstyle{plain}
\newtheorem{thm}{Theorem}[section]

\theoremstyle{definition}
\newtheorem{dfn}{Definition}[section]
\newtheorem{rem}{Remark}[section]
\newtheorem{asm}{Assumption}[section]

\newcounter{num}

\newfont{\bg}{cmr9 scaled\magstep4}
\newcommand{\bigzerol}{\smash{\lower1.0ex\hbox{\bg 0}}}

\newcommand{\R}{\mathbb{R}}
\newcommand{\bracket}[1]{\langle #1 \rangle}

\begin{document}
	
\title[Optimal decay rate for Klein-Gordon systems]{Optimal decay rate of solutions for nonlinear Klein-Gordon systems of critical type}

\author[S. Masaki]{Satoshi MASAKI}
\address[]{Division of Mathematical Science, Department of Systems Innovation, Graduate School of Engineering Science, Osaka University, Toyonaka, Osaka, 560-8531, Japan}
\email{masaki@sigmath.es.osaka-u.ac.jp}

\author[K. Sugiyama]{Koki SUGIYAMA}
\address[]{Division of Mathematical Science, Department of Systems Innovation, Graduate School of Engineering Science, Osaka University, Toyonaka, Osaka, 560-8531, Japan}
\email{k-sugiyama@sigmath.es.osaka-u.ac.jp}

\date{}
	
\maketitle

\begin{abstract}
	We consider the decay rate of solutions to nonlinear Klein-Gordon systems with a critical type nonlinearity. We will specify the optimal decay rate for a specific class of Klein-Gordon systems containing the dissipative nonlinearites. It will turn out that the decay rate which is previously found in some models is optimal.
\end{abstract}

\section{Introduction}
In this article,
we consider the following system of nonlinear Klein-Gordon equations:
\begin{equation}
 (\Box +m_j^2)u_j = F_j(u,\partial_t u,Du), \ \ (t,x)\in (0,\infty)\times \mathbb{R}^d, \  j=1,\cdots, N, \label{eqn} 
\end{equation}
where $\Box = \partial_t^2 - \Delta$, $d$ is a positive integer,  $u=(u_j)_{1\leq j\leq N}$ is an $\mathbb{R}^N$-valued unknown function, and
\begin{equation*}
 Du:=(\partial_i u_j)_{\substack{1\leq i\leq d \\1\leq j\leq N}}
\end{equation*}
with $\partial_k=\partial_{x_k}(k=1,\cdots,d)$. The masses $m_j$ are positive constants. We assume that the nonlinearity $F=(F_j)_{1\leq j \leq N}$ is of critical order, that is, it satisfies 
\begin{equation}
 |F(\zeta, \eta, \Theta)| \leq C\left( |\zeta|^2+|\eta|^2 +| \Theta |_{Fr}^2 \right)^{\frac12(1+\frac{2}{d})} \label{assume_nt}
\end{equation}
for $(\zeta, \eta, \Theta) \in \mathbb{R}^N\times \mathbb{R}^N\times \mathbb{R}^{N\times d}$.
Here $|\cdot |_{Fr}$ denotes the Frobenius norm,
\[
 | \Theta |_{Fr}^2 = \sum_{i=1}^{N} \sum_{j=1}^d |\theta_{ij}|^2
\]
of a matrix $\Theta=(\theta_{ij})_{1\leq i\leq N, 1\leq j\leq d}$.
We note that \eqref{eqn} includes a single Klein-Gordon equation as the special case $N=1$.
Here we are interested in the decay rate of solutions to \eqref{eqn}.
Our main goal is to give an upper bound on the decay rate of solutions to a class of Klein-Gordon equations and/or systems.
In particular, it will turn out that, in some specific models, a known decay rate is an optimal one.

There is a number of previous results on the decay rate and asymptotic behavior of Klein-Gordon equations and systems.
It is known that the order of the nonlinearity in \eqref{assume_nt} is a critical order regarding the long-time asymptotic behavior.

Let us introduce some previous results in the single equation case.
First, a solution of a single linear Klein-Gordon equation decays like
\begin{equation}
	\|u(t)\|_{L^\infty} = O(t^{-\frac{d}2}). \label{decay_free}
\end{equation}
As for the single equation with the power type nonlinearity $F=|u|^pu$, there exists a solution which approaches to a free solution as $t\rightarrow \infty$ if $p>2/d$ (see \cite{Klainerman-Ponce1983,Klainerman1982,Shatah1982,Shatah1985}).
Remark that the asymptotics holds in $L^\infty$ in some cases and that
the decay rate of a solution is the same as (\ref{decay_free}) in such cases.
On the other hand, if $0<p \leq 2/d$, Glassey \cite{Glassey1973} and Matsumura \cite{Matsumura1976} show that there is no non-trivial solution which asymptotically approach to a free solution. 
Our case $p=2/d$ is critical in this sense.

When the nonlinearity is of critical order, there are several possibilities on the asymptotic behavior of solutions and it depends on the shape of the nonlinearity.
For instance, as for the equation 
\begin{eqnarray}
(\Box + 1)u = F(u, \partial_t u, \partial_x u), \  \  (t,x)\in (0,\infty)\times \mathbb{R}, \nonumber
\end{eqnarray}
which is the case of $d=1$ and $N=1$ in \eqref{eqn}, the following is known.
Georgiev-Yordanov \cite{Georgiev-Yordanov1996} consider the cubic nonlinearity $F=u^3$ and prove that there is a global solution which decays like $O(t^{-\frac12})$ in $L^{\infty}$ but does not behave like a free solution.
In this case, it is known that the asymptotic behavior of the solution is a modified scattering (see \cite{Delort2001,Hayashi-Naumkin2008}).
On the other hand, Moriyama \cite{Moriyama1997} and Katayama \cite{Katayama1999} shows that there are some equations such that a solution behaves like a free solution (for example, a equation with $F=3u u_t^2 - 3u u_x^2 - u^3$).

Not only the asymptotic behavior but also the decay rate is different from (\ref{decay_free}) in some cases.
Indeed, if we take $F=-(\partial_t u)^3$, a logalithmic decay 
\[
	\|u(t)\|_{L^\infty} = O(t^{-\frac12}(\log t)^{-\frac12})
\]
is obtained by Sunagawa \cite{Sunagawa2006}.
Note that the decay rate of the solution is faster than a free solution due to the influence of the nonlinearity.
This kind of phenomenon is due to so-called nonlinear dissipation.

Next we consider the system case. In general, the behavior of the solution is expected to become richer in the system case.
The asymptotic behavior in the single case described above is also found in the system case.
First, Sunagawa \cite{Sunagawa2003} and Katayama-Ozawa-Sunagawa \cite{Katayama-Ozawa-Sunagawa2012} find some conditions on the nonlinearity $F$ and mass $m_j$ which ensure the existence of an asymptotically free solution for $d=1,2$.
At the same time, Sunagawa \cite{Sunagawa2004} shows that there exists a solution which decays like $O(t^{-\frac12})$ in $L^\infty$, the same decay as a free solution, but does not behave like a free solution and exhibits a modified-scattering type behavior.
Masaki-Segata-Uriya \cite{Masaki-Segata-Uriya2018} show modified scattering in the complex-valued case for $d=1,2$.
Remark that the single complex-valued equation is equivalent to a system of real-valued system (see also \cite{Sunagawa2005-2}).
In the above cases, the decay rate of the solutions is the same as that of \eqref{decay_free}.

Kim-Sunagawa \cite{Kim-Sunagawa2014} show that if we consider the nonlinearity
\begin{equation}\label{E:example}
	F_j=\mu_1 |u|^2 u_j-\mu_2 |\partial_t u|^2\partial_t u_j \  \  (\mu_1 \in \mathbb{R},\, \mu_2>0,\, j = 1,2,\dots,N),
\end{equation}
then the corresponding model
admits a solution which is compactly-supported and decays like
\begin{eqnarray}
	\| u(t) \|_{L^\infty} = O(t^{-\frac12}(\log{t})^{-\frac12}) \label{decay}
\end{eqnarray}
under an appropriate condition on the initial data.
Moreover, the solution $u$ satisfies 
\begin{equation}
	\| u(t,\cdot) \|_{L^\infty} + \| \partial_t u(t,\cdot) \|_{L^\infty} + \| \partial_x u(t,\cdot) \|_{L^\infty} \leq Ct^{-\frac{1}{2}}(\log{t})^{-\frac{1}{2}} \label{dissipation_est}
\end{equation}
in this case.
We note that this nonlinearity is a typical example which has the nonlinear dissipation effect.
They give a condition on nonlinearity for which there exists a solution satisfying \eqref{dissipation_est}.

The ratio of the mass coefficients $m_j$ also matters in the system case. \cite{Sunagawa2003} gave such an example. It is shown that,
for one fixed nonlinearity, the corresponding solution blows up in finite time for some choices of $m_j$, and the global solution exists for other choices.
According to \cite{Sunagawa2005}, the system
\[
	\left\{
	\begin{array}{l}
		(\Box + m_1^2)u_1=0 \\
		(\Box + m_2^2)u_2 =u_1^3
	\end{array}
	\right.
\]
is an example for which the solution behaves differently from a solution of a single equation.
The second component $u_2$ of a solution decays no faster than $t^{-\frac12}(\log{t})^{\frac12}$ in $L^{\infty}$.

In this article, we show that a class of the nonlinear Klein-Gordon equations/system can not admit a (good) solution decaying faster than the decay rate (\ref{dissipation_est}).
This in particular shows that the rate \eqref{dissipation_est} is optimal for a class of nonlinearity including \eqref{E:example}.
Recently, Kita \cite{Kita2019} studies the optimal decay rate of solution to a dissipative nonlinear Schr\"odinger equation.
Our approach is inspired by his argument.


\subsection{Main results}

Let us state our main theorem.
The assumptions on a nonlinearity is as follows:

\begin{asm} \label{assume1}
The nonlinearity $F$ satisfies the estimate (\ref{assume_nt}).
Moreover, we have uniqueness of 
a classical solution to the initial value problem of the equation \eqref{eqn} with $F$, that is, 
if $u_1$ and $u_2$ are two classical solutions and if $u_1(t_0)=u_2(t_0)$ holds for some $t_0$ then we have $u_1 \equiv u_2$.
\end{asm}

We say a solution is forward-global if the solution exists on $[T,\infty)$ for some $T\in \mathbb{R}$.
Now we state our main theorem.

\begin{thm} \label{thm1}
Suppose that the nonlinearity $F$ satisfies Assumption  \ref{assume1}.
	Let $u$ be a forward-global classical solution to (\ref{eqn}) with the nonlinearity $F$. Assume that there exist constants $C>0$ and $T_1>0$ such that 
	\begin{eqnarray}
		\| u(t,\cdot) \|_{L^\infty} + \| |Du(t,\cdot)|_{Fr} \|_{L^\infty} \leq Ct^{-\frac{d}{2}}(\log{t})^{-\frac{d}{2}} \label{assume_other}
	\end{eqnarray}
	holds for all $t\geq T_1$. Also assume that for any $\varepsilon  > 0$ there exists $T_2 \geq T_1$ such that 
	\begin{eqnarray}
		\| \partial_{t} u(t,\cdot) \|_{L^\infty} \leq \varepsilon t^{-\frac{d}{2}}(\log{t})^{-\frac{d}{2}} \label{assume_ut}
	\end{eqnarray}
	holds for all $t\geq T_2$. Then $u$ is identically zero.
\end{thm}

\begin{rem}
	In Theorem \ref{thm1}, we assume that the existence of a forward-global solution of (\ref{eqn}). We note that the existence of a forward-global solution is not trivial. 
	In fact, even in the simple case of $d=1$ and $N=1$, small data global existence is not obvious (for instance, see \cite{Katayama1999}).
	A sufficient condition for small data global existence is given in \cite{Delort2001}.
\end{rem}

\begin{rem}
	We would emphasize that
	it is not necessary to assume that the solution $u$ is compactly supported.
\end{rem}

We consider a classical solution in Theorem \ref{thm1}.
However, by combining with an appropriate well-posedness result, we can show the same conclusion as Theorem \ref{thm1} for other class of solutions.
As such an example, we give a result for $H_x^2\times H_x^1$-solutions, which is our second result.
In the following, we assume $d=1$ and $m_1=\cdots = m_N$ for simplicity.

Before stating the result, let us introduce the notion of an $H^2_x \times H^1_x$-solution.

\begin{dfn} \label{dfn_sol}
	We say that a function $u$ is an $H_x^2\times H_x^1$-solution on an interval $I\subset \R$
to \eqref{eqn} if $u \in C_t(I;(H_x^2(\mathbb{R}))^N)\cap C_t^1(I;(H_x^1(\mathbb{R}))^N)$ and u obeys the Duhamel formula:
	\begin{eqnarray}
		\left( \begin{array}{c}
		u_j(t) \\
		\partial_t u_j(t)
		\end{array} \right) = \left( \begin{array}{cc}
		\cos{\bracket{\partial_x} (t-t_0)} & \bracket{\partial_x}^{-1} \sin{\bracket{\partial_x}(t-t_0)} \\
		-\bracket{\partial_x}\sin{\bracket{\partial_x}(t-t_0)} & \cos{\bracket{\partial_x}(t-t_0)}
		\end{array}	\right) \left( \begin{array}{c}
		u_j(t_0) \\
		\partial_t u_j(t_0)
		\end{array} \right) \nonumber \\
		+ \int_{t_0}^{t} \left( \begin{array}{c}
		\bracket{\partial_x}^{-1}\sin{\bracket{\partial_x}(t-s)} \\
		\cos{\bracket{\partial_x}(t-s)}
		\end{array} \right) F_j(u(s), \partial_tu(s), \partial_x u(s)) ds \label{duhamel}
	\end{eqnarray}
	for $t_0 ,t\in I$.
\end{dfn}

We state the assumption on nonlinearities.

\begin{asm} \label{assume2}
	The nonlinearity $F$ satisfies
	\begin{eqnarray}
	|\partial_{\zeta, \eta, \theta}^\alpha F(\zeta, \eta, \theta)| \leq C \left( |\zeta|^2 + |\eta|^2 + |\theta|^2 \right)^{\frac{3-|\alpha|}{2}} \label{assume_wp}
	\end{eqnarray}
	for $0\leq |\alpha| \leq 3$.
\end{asm}

By a standard argument, one can show that if the nonlinearity $F$ satisfies the above assumption then
the local well-posedness of the initial value problem of \eqref{eqn} holds in $H_x^2 \times H_x^1$.
Furthermore, the equation admits a classical solution for smooth data.

The following is our second result.

\begin{thm} \label{thm2}
	Suppose that the nonlinearity $F$ satisfies Assumption \ref{assume2}.
	Let $u$ be a forward-global $H_x^2\times H_x^1$-solution to \eqref{eqn}. If $u$ satisfies \eqref{assume_other} and \eqref{assume_ut}
	then $u$ is identically zero.
\end{thm}

It is clear that
Kim-Sunagawa's model \eqref{E:example} satisfies Assumption \ref{assume2}. 
The above theorem shows that the decay rate \eqref{dissipation_est} is optimal for this model.


\section{Proof of Theorem \ref{thm1}}

In the next section, we prove Theorem \ref{thm1}. 
A key ingredient is a \emph{localized linear energy}: For a classical solution $u$,
$t\in \mathbb{R}$, and $y \in \mathbb{R}^d$,
	\begin{equation}\label{E:lle}
	E_{t,y}\left(u\left( t \right)\right) := \sum_{j=1}^N \int_{|x-y|\leq t} \left\{ |\nabla  u_j(t,x)|^2+|\partial_t u_j(t,x)|^2+m_j^2 |u_j(t,x)|^2 \right\}dx.
	\end{equation}
Remark that the quantity is well-defined for a classical solution even when it does not decay near the spatial infinity.

\begin{proof}
	Define the localized linear energy $E_{t,y}$
	for any $y\in \mathbb{R}$ as in \eqref{E:lle}. Then, it follows that
	\begin{eqnarray*}
	&&\frac{d}{dt}E_{t,y}(u(t)) = 2\sum_{j=1}^N \int_{|x-y|\leq t}\left\{ (\nabla u_j)\cdot (\partial_t\nabla u_j)+(\partial_t u_j)(\partial_t^2 u_j) + (m_j^2 u_j) (\partial_t u_j) \right\}dx \label{first_term}\\ && \hspace{3cm} +\sum_{j=1}^N \int_{|x-y|=t}\left\{ |\nabla  u_j(t,x)|^2+|\partial_t u_j(t,x)|^2+m_j^2 |u_j(t,x)|^2 \right\}dx.
	\end{eqnarray*}
	We calculate the first term above. By (\ref{assume_nt}), we deduce \begin{eqnarray*}
	&&\int_{|x-y|\leq t}\left\{ (\nabla u_j)\cdot(\partial_t\nabla u_j)+(\partial_t u_j)(\partial_t^2 u_j) + (m_j^2 u_j)(\partial_t u_j) \right\} dx \\
	&=& \int_{|x-y|\leq t}(\partial_t u_j) \left( \partial_t^2-\Delta + m_j^2 \right)u_j dx + \int_{|x-y|=t} \partial_t u_j \nabla u_j \cdot \frac{x}{|x|} dx \\
	&=& \int_{|x-y|\leq t}(\partial_t u_j) F_j(u,\partial_t u, Du)dx + \int_{|x-y|=t} \partial_t u_j \nabla u_j \cdot \frac{x}{|x|} dx \\
	&\geq& -C \int_{|x-y|\leq t}(\partial_t u_j) \left( |u|^2+|\partial_t u|^2+|D u|_{Fr}^2\right)^{\frac12(1+\frac{2}{d})} dx \label{calc_first_term} \\ && \hspace{5cm} + \int_{|x-y|=t} \partial_t u_j \nabla u_j \cdot \frac{x}{|x|} dx .
	\end{eqnarray*}
	Using the assumption (\ref{assume_other}) and (\ref{assume_ut}), we have
	\begin{eqnarray*}
	&& -C \int_{|x-y|\leq t}(\partial_t u_j)\left( |u|^2+|\partial_t u|^2+|D u|_{Fr}^2\right)^{\frac12(1+\frac{2}{d})} dx \\
	&\geq& -C \| \partial_t u \|_{L^\infty} \left(  \| u \|_{L^\infty} + \| \partial_t u \|_{L^\infty} + \| |D u|_{Fr} \|_{L^\infty} \right)^{-1+\frac{2}{d}} E_{t,y}(u(t)) \\
	&\geq& -C \varepsilon t^{-1} (\log{t})^{-1} E_{t,y}(u(t))
	\end{eqnarray*}
	when $d=1,2$. Similarly, one has
	\begin{eqnarray*}
		&& -C \int_{|x-y|\leq t}(\partial_t u_j)\left( |u|^2+|\partial_t u|^2+|D u|_{Fr}^2\right)^{\frac12(1+\frac{2}{d})} dx \\
		&\geq& -C \| \partial_t u \|_{L^\infty}^{\frac2d} \int_{|x-y|\leq t}|\partial_t u|^{1-\frac2d}\left( |u|^2+|\partial_t u|^2+|D u|_{Fr}^2\right)^{\frac12(1+\frac{2}{d})} dx \\
		&\geq& -C \varepsilon^{\frac2d} t^{-1} (\log{t})^{-1} E_{t,y}(u(t)) 
	\end{eqnarray*}
	when $d \geq 3$. 
	Combining these estimates, we obtain
	\begin{eqnarray*}
	\frac{d}{dt}E_{t,y}(u(t)) &\geq& -C\varepsilon t^{-1} (\log{t})^{-1} E_{t,y}(u(t)) \\ &&\hspace{3mm} +\sum_{j=1}^N \int_{|x-y|=t}\left\{m_j^2 |u_j|^2 + |\nabla  u_j|^2+|\partial_t u_j|^2+2\partial_t u_j \nabla u_j \cdot \frac{x}{|x|} \right\}dx \\
	&\geq& -C\varepsilon t^{-1} (\log{t})^{-1} E_{t,y}(u(t)) \\ &&\hspace{3mm} +\sum_{j=1}^N \int_{|x-y|=t}\left\{m_j^2 |u_j|^2 + \left| \partial_t u_j \frac{x}{|x|}+\nabla u_j \right|^2 \right\}dx.
	\end{eqnarray*}
	Therefore we reach to the estimate
	\begin{eqnarray}
	\frac{d}{dt}E_{t,y}(u(t)) \geq -C\varepsilon t^{-1} (\log{t})^{-1} E_{t,y}(u(t)). \label{energy_est}
	\end{eqnarray}
	Fix $\delta \in (0,d)$.
	By using (\ref{energy_est}), we deduce that
	\begin{eqnarray*}
	\frac{d}{dt}\biggl( (\log{t})^\delta E_{t,y}(u(t)) \biggr) &=& \delta t^{-1}(\log{t})^{\delta-1}E_{t,y}(u(t)) + (\log{t})^\delta \frac{d}{dt}E_{t,y}(u(t)) \\
	&\geq& \delta t^{-1}(\log{t})^{\delta-1}E_{t,y}(u(t)) -C\varepsilon t^{-1} (\log{t})^{\delta-1} E_{t,y}(u(t)) \\
	&\geq& 0
	\end{eqnarray*}
	for $\varepsilon$ sufficiently small and $t \geq T_2(\varepsilon)$. Integrating this with respect to time, we have
	\begin{equation*}
	(\log{t_1})^\delta E_{t_1,y}(u(t_1)) \leq (\log{t_2})^\delta E_{t_2,y}(u(t_2))
	\end{equation*}
	for $T_2 \leq t_1 \leq t_2$. 

	On the other hand, by using (\ref{assume_other}) and (\ref{assume_ut}), we deduce 
	\begin{eqnarray*}
	E_{t,y}(u(t)) &\leq& C \left( \| u \|_{L^\infty} + \| \partial_t u \|_{L^\infty} + \| |Du|_{Fr} \|_{L^\infty} \right)^{2} \sum_{j=1}^N \int_{|x-y|\leq t}dx \\
	&\leq & C  (\log{t})^{-d}
	\end{eqnarray*}
	for $t\ge T_2$.
	Thus, we obtain
	\begin{eqnarray*}
	(\log{t_1})^\delta E_{t_1,y}(u(t_1)) &\leq& C (\log{t_2})^{\delta-d} .
	\end{eqnarray*}
	Since $\delta <d$,
	the right-hand side converges to zero as $t_2$ tends to infinity, showing that $E_{t_1,y}(u(t_1)) = 0$.
	One then sees that
	\begin{eqnarray*}
	u(t_1,x)=0, \  \  \partial_t u(t_1,x) = 0
	\end{eqnarray*}
	for $|x-y|<t_1$. Since $t_1$ is independent on $y$, this is true for any choice of $y \in \mathbb{R}^d$.
	Thus, we obtain
	\begin{eqnarray}
	u(t_1,x)=0, \  \  \partial_t u(t_1,x) = 0 \label{conclusion}
	\end{eqnarray}
	for all $x\in\mathbb{R}^d$. By the assumption \eqref{assume_nt}, the zero solution is also a forward-global classical solutions
 satisfying (\ref{conclusion}). Hence, we obtain the desired conclusion by the uniqueness of a classical solution.
\end{proof}

\begin{section}{Proof of Theorem \ref{thm2}}
	
	\begin{proof}
		We set the localized linear energy $E_{t,y}$ as in \eqref{E:lle}.
		Pick sequences $\{u_{0,n}\}\subset H_x^3(\R)^N$ and $\{u_{1,n}\}\subset H_x^2(\R)^N$ so that
		\[
			\| u(0) - u_{0,n}(0) \|_{H_x^2} + \| \partial_t u(0) - \partial_t u_{1,n}(0) \|_{H_x^1}  \rightarrow 0 
		\]
		as $n\rightarrow \infty$. 
	By Assumption \ref{assume2}, there exists a classical solution $u_n(t)$ to \eqref{eqn} with data $(u_n(0), \partial u_n(0)) = (u_{0,n}, u_{1,n})$.
		Then, a standard blowup criterion and continuous dependence show that 
		\begin{equation}
			\| u_n - u \|_{L_t^{\infty}((0,\tau);H_x^2)} + \| \partial_t u_n - \partial_t u \|_{L_t^{\infty}((0,\tau);H_x^1)} \rightarrow 0 \label{continuity}
		\end{equation}
		as $n\rightarrow \infty$ for all $\tau > 0$. 
		
		Let $\varepsilon>0$ to be chosen later.
		Fix $\tau > 2T_2(\varepsilon)$.
		By \eqref{continuity}, we deduce that there exists a constant $N_0=N_0(\varepsilon,\tau)\in \mathbb{N}_{\geq 0}$ such that
		\begin{eqnarray*}
			\| \partial_t u_n(t) \|_{L_x^\infty} &\leq& \| \partial_t u_n(t)-\partial_t u(t) \|_{L_x^\infty} +\| \partial_t u(t) \|_{L_x^\infty}\\
			&\leq& C \| \partial_t u_n(t)-\partial_t u(t) \|_{H_x^1} +\| \partial_t u(t) \|_{L_x^\infty}\\
			&\leq& C\varepsilon t^{-\frac12}(\log{t})^{-\frac12}
		\end{eqnarray*}
		for all $n>N_0$ and $t\in(2T_2, \tau)$. In the same way, we also deduce that there exists $N_1\ge N_0$ such that
		\begin{eqnarray*}
			\| u_n(t) \|_{L_x^\infty} + \| \partial_x u_n(t) \|_{L_x^\infty} &\leq& C \| u_n(t)- u(t) \|_{H_x^2} +\| \partial_x u(t) \|_{L_x^\infty} +\| u(t) \|_{L_x^\infty} \\
 			&\leq& C t^{-\frac12}(\log{t})^{-\frac12}
		\end{eqnarray*}
		for all $n>N_1$ and $t\in(2T_2, \tau)$. 
		
		The rest of the proof is almost the same as that of Theorem \ref{thm1}.
		By arguments similar to the proof of \eqref{energy_est}, we obtain 
		\begin{eqnarray*}
			\frac{d}{dt}E_{t,y}(u_n(t)) \geq -C\varepsilon t^{-1} (\log{t})^{-1} E_{t,y}(u_n(t))
		\end{eqnarray*}
		for $t\in(2T_2, \tau)$ and $n>N_1$.	Similarly,
		\begin{eqnarray*}
			E_{t,y}(u_n(t)) \leq C(\log{t})^{-d}
		\end{eqnarray*}
		holds for $t\in(2T_2, \tau)$ and $n>N_1$. 
		We now fix $\delta \in (0,d)$ and choose $\varepsilon>0$ so small that
		we have
		\begin{eqnarray*}
			(\log{2T_2})^\delta E_{2T_2,y}(u_n(2T_2)) &\leq& (\log{\tau})^\delta E_{\tau,y}(u_n(\tau)) \\ &\leq&(\log{\tau})^{\delta-d}
		\end{eqnarray*}
		for all $n>N_1(\varepsilon,\tau)$.
		We pass to the limit $n\to\infty$ to obtain
		\[
			(\log{2T_2})^\delta E_{2T_2,y}(u(2T_2)) \leq (\log{\tau})^{\delta-d}.
		\]
		Since $\tau > 2T(\varepsilon)$ is arbitrary, we obtain
		\[
			(\log{2T_2})^\delta E_{2T_2,y}(u(2T_2)) = 0
		\]
		by letting $\tau \rightarrow \infty$. Since this is true for all $y \in \mathbb{R}^d$, we have
		$u(2T_2)=\partial_t u(2T_2)=0$. By the uniqueness of an $H^2_x\times H^1_x$-solution, we conclude that
		$u$ is identically zero.
	\end{proof}
\end{section}

\subsection*{Acknowledgment}
S.M. was supported by JSPS KAKENHI Grant Numbers JP17K14219, JP17H02854, JP17H02851, and JP18KK0386.



\begin{bibdiv}
\begin{biblist}

\bib{Delort2001}{article}{
      author={Delort, Jean-Marc},
       title={Existence globale et comportement asymptotique pour
  l'\'{e}quation de {K}lein-{G}ordon quasi lin\'{e}aire \`a donn\'{e}es petites
  en dimension 1},
        date={2001},
        ISSN={0012-9593},
     journal={Ann. Sci. \'{E}cole Norm. Sup. (4)},
      volume={34},
      number={1},
       pages={1\ndash 61},
  url={https://doi-org.remote.library.osaka-u.ac.jp:8443/10.1016/S0012-9593(00)01059-4},
      review={\MR{1833089}},
}

\bib{Georgiev-Yordanov1996}{article}{
      author={Georgiev, Vladimir},
      author={Yordanov, Borislav},
       title={Asymptotic behaviour of the one-dimensional klein--gordon
  equation with a cubic nonlinearity},
        date={1996},
     journal={preprint},
}

\bib{Glassey1973}{article}{
      author={Glassey, Robert~T.},
       title={On the asymptotic behavior of nonlinear wave equations},
        date={1973},
        ISSN={0002-9947},
     journal={Trans. Amer. Math. Soc.},
      volume={182},
       pages={187\ndash 200},
  url={https://doi-org.remote.library.osaka-u.ac.jp:8443/10.2307/1996529},
      review={\MR{330782}},
}

\bib{Hayashi-Naumkin2008}{article}{
      author={Hayashi, Nakao},
      author={Naumkin, Pavel~I.},
       title={The initial value problem for the cubic nonlinear
  {K}lein-{G}ordon equation},
        date={2008},
        ISSN={0044-2275},
     journal={Z. Angew. Math. Phys.},
      volume={59},
      number={6},
       pages={1002\ndash 1028},
  url={https://doi-org.remote.library.osaka-u.ac.jp:8443/10.1007/s00033-007-7008-8},
      review={\MR{2457221}},
}

\bib{Katayama1999}{article}{
      author={Katayama, Soichiro},
       title={A note on global existence of solutions to nonlinear
  {K}lein-{G}ordon equations in one space dimension},
        date={1999},
        ISSN={0023-608X},
     journal={J. Math. Kyoto Univ.},
      volume={39},
      number={2},
       pages={203\ndash 213},
  url={https://doi-org.remote.library.osaka-u.ac.jp:8443/10.1215/kjm/1250517908},
      review={\MR{1709289}},
}

\bib{Katayama-Ozawa-Sunagawa2012}{article}{
      author={Katayama, Soichiro},
      author={Ozawa, Tohru},
      author={Sunagawa, Hideaki},
       title={A note on the null condition for quadratic nonlinear
  {K}lein-{G}ordon systems in two space dimensions},
        date={2012},
        ISSN={0010-3640},
     journal={Comm. Pure Appl. Math.},
      volume={65},
      number={9},
       pages={1285\ndash 1302},
  url={https://doi-org.remote.library.osaka-u.ac.jp:8443/10.1002/cpa.21392},
      review={\MR{2954616}},
}

\bib{Kim-Sunagawa2014}{article}{
      author={Kim, Donghyun},
      author={Sunagawa, Hideaki},
       title={Remarks on decay of small solutions to systems of
  {K}lein-{G}ordon equations with dissipative nonlinearities},
        date={2014},
        ISSN={0362-546X},
     journal={Nonlinear Anal.},
      volume={97},
       pages={94\ndash 105},
  url={https://doi-org.remote.library.osaka-u.ac.jp:8443/10.1016/j.na.2013.11.025},
      review={\MR{3146374}},
}

\bib{Kita2019}{article}{
      author={Kita, Naoyasu},
       title={Optimal decay rate of solutions to 1d {S}chr\"{o}dinger equation
  with cubic dissipative nonlinearity},
        date={2019},
     journal={Journal of Applied Science and Engineering A},
      volume={1},
      number={1},
       pages={15\ndash 18},
}

\bib{Klainerman-Ponce1983}{article}{
      author={Klainerman, S.},
      author={Ponce, Gustavo},
       title={Global, small amplitude solutions to nonlinear evolution
  equations},
        date={1983},
        ISSN={0010-3640},
     journal={Comm. Pure Appl. Math.},
      volume={36},
      number={1},
       pages={133\ndash 141},
  url={https://doi-org.remote.library.osaka-u.ac.jp:8443/10.1002/cpa.3160360106},
      review={\MR{680085}},
}

\bib{Klainerman1982}{article}{
      author={Klainerman, Sergiu},
       title={Long-time behavior of solutions to nonlinear evolution
  equations},
        date={1982},
        ISSN={0003-9527},
     journal={Arch. Rational Mech. Anal.},
      volume={78},
      number={1},
       pages={73\ndash 98},
  url={https://doi-org.remote.library.osaka-u.ac.jp:8443/10.1007/BF00253225},
      review={\MR{654553}},
}

\bib{Masaki-Segata-Uriya2018}{article}{
      author={Masaki, Satoshi},
      author={Segata, Junichi},
      author={Uriya, Kota},
       title={Long range scattering for the complex-valued klein-gordon
  equation with quadratic nonlinearity in two dimensions},
        date={2018},
     journal={arXiv preprint arXiv:1810.02158},
}

\bib{Matsumura1976}{article}{
      author={Matsumura, Akitaka},
       title={On the asymptotic behavior of solutions of semi-linear wave
  equations},
        date={1976/77},
        ISSN={0034-5318},
     journal={Publ. Res. Inst. Math. Sci.},
      volume={12},
      number={1},
       pages={169\ndash 189},
  url={https://doi-org.remote.library.osaka-u.ac.jp:8443/10.2977/prims/1195190962},
      review={\MR{0420031}},
}

\bib{Moriyama1997}{article}{
      author={Moriyama, Kazunori},
       title={Normal forms and global existence of solutions to a class of
  cubic nonlinear {K}lein-{G}ordon equations in one space dimension},
        date={1997},
        ISSN={0893-4983},
     journal={Differential Integral Equations},
      volume={10},
      number={3},
       pages={499\ndash 520},
      review={\MR{1744859}},
}

\bib{Shatah1982}{article}{
      author={Shatah, Jalal},
       title={Global existence of small solutions to nonlinear evolution
  equations},
        date={1982},
        ISSN={0022-0396},
     journal={J. Differential Equations},
      volume={46},
      number={3},
       pages={409\ndash 425},
  url={https://doi-org.remote.library.osaka-u.ac.jp:8443/10.1016/0022-0396(82)90102-4},
      review={\MR{681231}},
}

\bib{Shatah1985}{article}{
      author={Shatah, Jalal},
       title={Normal forms and quadratic nonlinear {K}lein-{G}ordon equations},
        date={1985},
        ISSN={0010-3640},
     journal={Comm. Pure Appl. Math.},
      volume={38},
      number={5},
       pages={685\ndash 696},
  url={https://doi-org.remote.library.osaka-u.ac.jp:8443/10.1002/cpa.3160380516},
      review={\MR{803256}},
}

\bib{Sunagawa2003}{article}{
      author={Sunagawa, Hideaki},
       title={On global small amplitude solutions to systems of cubic nonlinear
  {K}lein-{G}ordon equations with different mass terms in one space dimension},
        date={2003},
        ISSN={0022-0396},
     journal={J. Differential Equations},
      volume={192},
      number={2},
       pages={308\ndash 325},
  url={https://doi-org.remote.library.osaka-u.ac.jp:8443/10.1016/S0022-0396(03)00125-6},
      review={\MR{1990843}},
}

\bib{Sunagawa2004}{article}{
      author={Sunagawa, Hideaki},
       title={A note on the large time asymptotics for a system of
  {K}lein-{G}ordon equations},
        date={2004},
        ISSN={0385-4035},
     journal={Hokkaido Math. J.},
      volume={33},
      number={2},
       pages={457\ndash 472},
  url={https://doi-org.remote.library.osaka-u.ac.jp:8443/10.14492/hokmj/1285766177},
      review={\MR{2073010}},
}

\bib{Sunagawa2005}{article}{
      author={Sunagawa, Hideaki},
       title={Large time asymptotics of solutions to nonlinear {K}lein-{G}ordon
  systems},
        date={2005},
        ISSN={0030-6126},
     journal={Osaka J. Math.},
      volume={42},
      number={1},
       pages={65\ndash 83},
  url={http://projecteuclid.org.remote.library.osaka-u.ac.jp/euclid.ojm/1153494315},
      review={\MR{2130963}},
}

\bib{Sunagawa2005-2}{article}{
      author={Sunagawa, Hideaki},
       title={Remarks on the asymptotic behavior of the cubic nonlinear
  {K}lein-{G}ordon equations in one space dimension},
        date={2005},
        ISSN={0893-4983},
     journal={Differential Integral Equations},
      volume={18},
      number={5},
       pages={481\ndash 494},
      review={\MR{2136975}},
}

\bib{Sunagawa2006}{article}{
      author={Sunagawa, Hideaki},
       title={Large time behavior of solutions to the {K}lein-{G}ordon equation
  with nonlinear dissipative terms},
        date={2006},
        ISSN={0025-5645},
     journal={J. Math. Soc. Japan},
      volume={58},
      number={2},
       pages={379\ndash 400},
  url={http://projecteuclid.org.remote.library.osaka-u.ac.jp/euclid.jmsj/1149166781},
      review={\MR{2228565}},
}

\end{biblist}
\end{bibdiv}

\end{document}